# A PECULIAR TWO POINT BOUNDARY VALUE PROBLEM

By Huadong Pang[1,2] and Daniel W. Stroock[2]

*Massachusetts Institute of Technology*

In this paper we consider a one-dimensional diffusion equation on the interval $[0, 1]$ satisfying non-Feller boundary conditions. As a consequence, the initial value Cauchy problem fails to preserve nonnegativity or boundedness. Nonetheless, probability theory plays an interesting role in our analysis and understanding of solutions to this equation.

**1. Introduction.** In this article, we continue the study, started in [6] and [7], of a diffusion equation in one dimension with a boundary condition for which the minimum principle fails. The main distinction between the situation here and the one studied earlier is that we are now dealing with a problem in which there are two boundary points, not just one, and the addition of the second boundary point introduces some new phenomena which we find interesting.

Although the relationship is not immediately apparent, related considerations appear in [3] and [4].

1.1. *The problem and a basic result.* Let F be the space of bounded functions on $[0, 1]$ which are continuous on $(0, 1)$ but not necessarily continuous at the boundary $\{0, 1\}$. Convergence of $\{f_n\}_1^\infty \subseteq F$ to $f$ in F means that $\{\|f_n\|_u\}_1^\infty$ is bounded, $f_n(x) \longrightarrow f(x)$ for each $x \in [0, 1]$ and uniformly for $x$ in compact subsets of $(0, 1)$.

In the next definition, and hereafter, we use the probabilistic convention of writing $u(t, x)$ where analysts would use $u(x, t)$. As usual,

$$\dot{u} \equiv \frac{\partial u}{\partial t}, \qquad u' \equiv \frac{\partial u}{\partial x} \quad \text{and} \quad u'' \equiv \frac{\partial^2 u}{\partial x^2}.$$

Received April 2006; revised May 2006.
[1]Supported by the MIT Mathematics Department.
[2]Supported by the NSF Grant DMS-02-44991.
*AMS 2000 subject classifications.* 60J60, 35K15, 35K35.
*Key words and phrases.* One-dimensional diffusion, minimum principle.

Now let $U$ be the space of functions $u \in C^{1,2}((0,\infty) \times [0,1]; \mathbb{R})$ with the properties that $u$ is bounded on $(0,1] \times [0,1]$ and, for each $0 < T_1 < T_2 < \infty$, $\dot u$, $u'$ and $u''$ are bounded on $[T_1, T_2] \times [0,1]$. Note that we are insisting that $u$ be $C^{1,2}$ right up to, and including, the spacial boundary $(0,\infty) \times \{0,1\}$.

Because its proof is more easily understood after seeing the proofs of the other results in this article, we have put the derivation of the following basic existence and uniqueness statement into an [Appendix](#) at the end of this article.

THEOREM 1.1. *Let $(\mu, \sigma) \in \mathbb{R}^2$ be given.*

(i) *Suppose that $u \in U$ satisfies*

(1.1)
$$\dot u = \tfrac{1}{2} u'' + \mu u' \quad \text{on } (0,\infty) \times (0,1),$$
$$\dot u(t,0) = -\sigma u'(t,0) \quad \text{and} \quad \dot u(t,1) = \sigma u'(t,1) \quad \text{for } t \in (0,\infty).$$

*If, as $t \searrow 0$, $u(t,\cdot)$ converges uniformly on compact subsets of $(0,1)$, then both $u(t,0)$ and $u(t,1)$ converge as $t \searrow 0$, and so $u(t,\cdot)$ converges in $F$.*

(ii) *Given $f \in F$, there is a unique $u_f \in U$ which satisfies ([1.1](#)) and the initial condition that, as $t \searrow 0$, $u(t,\cdot)$ converges to $f$ in $F$.*

*In particular, if $Q_t f \equiv u_f(t, \cdot)$, then $\{Q_t : t \geq 0\}$ is a semigroup of bounded, continuous operators on $F$. [See ([3.2](#)) below for more information.]*

For semigroup enthusiasts, it may be helpful to think of the operator $Q_t$ as $\exp(t\mathcal{H})$ where $\mathcal{H} f = \tfrac{1}{2} f'' + \mu f'$ with domain

$\mathrm{dom}(\mathcal{H})$
$$= \{f \in C^2([0,1]; \mathbb{R}) : \tfrac{1}{2} f''(k) + \mu f'(k) = (-1)^{1-k} \sigma f'(x) \text{ for } k \in \{0,1\}\}.$$

For probabilists, it may be helpful to remark that, unless $\sigma \leq 0$, $\{Q_t : t \geq 0\}$ is *not* a Markov semigroup.

1.2. *Nonnegativity and growth of solutions.* If $\sigma \leq 0$, then $u_f(\cdot, \cdot) \geq 0$ if and only if $f \geq 0$, and therefore $\{Q_t : t \geq 0\}$ is a Markov (i.e., nonnegativity preserving) semigroup. This may be proved by either an elementary minimum principle argument or the well-known probabilistic model. [The corresponding diffusion is Brownian motion in $(0,1)$ with drift $\mu$ which, depending on whether $\sigma = 0$ or $\sigma < 0$, is either absorbed when it hits $\{0,1\}$ or has a "sticky" reflection there.] However, when $\sigma > 0$, the minimum principle is lost, and, as a consequence $\{Q_t : t \geq 0\}$ is no longer Markov. Nonetheless, we will show that there is a certain $\{Q_t : t \geq 0\}$-invariant subspace of $F$ on which the $Q_t$'s do preserve nonnegativity. To describe this subspace, we need the following.

THEOREM 1.2. *Given a continuously differentiable function $J:[0,1] \longrightarrow \mathbb{R}^2$, set*

$$B(J) = \begin{pmatrix} -2\mu\sigma + \frac{J_0'(0)}{2} & -\frac{J_0'(1)}{2} \\ \frac{J_1'(0)}{2} & 2\mu\sigma - \frac{J_1'(1)}{2} \end{pmatrix}$$

$$= \left(2\sigma\mu - \frac{1}{2}(J'(0), J'(1))\right)\begin{pmatrix} -1 & 0 \\ 0 & 1 \end{pmatrix}.$$

*Then, for each $\sigma > 0$ and $\mu \in \mathbb{R}$, there exist a unique solution $J^{\sigma,\mu}$ to*

(R)
$$\tfrac{1}{2}J''(x) - \mu J'(x) + B(J)J(x) = 0 \quad \text{on } [0,1]$$

$$J(0) = \begin{pmatrix} 2\sigma \\ 0 \end{pmatrix} \quad \text{and} \quad J(1) = \begin{pmatrix} 0 \\ 2\sigma \end{pmatrix}$$

*which satisfies*

(1.2) $$\max_{k \in \{0,1\}} \int_0^1 |J_k(x)|\, dx \begin{cases} \leq 1, & \text{if } \sigma \geq \mu \coth \mu, \\ < 1, & \text{if } \sigma < \mu \coth \mu. \end{cases}$$

*Moreover, $J^{\sigma,\mu} \geq 0$ in the sense that both of its components are nonnegative. Finally, set $B^{\sigma,\mu} = B(J^{\sigma,\mu})$. Then $B^{\sigma,\mu}$ has real eigenvalues $\lambda_1^{\sigma,\mu} < \lambda_0^{\sigma,\mu} \leq 0$, $\lambda_0^{\sigma,\mu} < 0$ if and only if $\sigma > \mu \coth \mu$, and the corresponding eigenvector $V_0^{\sigma,\mu}$ can be chosen to be strictly positive with $(V_0^{\sigma,\mu})_0 + (V_0^{\sigma,\mu})_1 = 1$, whereas the eigenvector $V_1^{\sigma,\mu}$ corresponding to $\lambda_1^{\sigma,\mu}$ can be chosen so that $(V_1^{\sigma,\mu})_0 > 0 > (V_1^{\sigma,\mu})_1$ and $(V_1^{\sigma,\mu})_0 - (V_1^{\sigma,\mu})_1 = 1$. (See Lemmas 2.1 and 2.2 below for more information.)*

Referring to the quantities in Theorem 1.2, we have the following. When $\mu = 0$, some of the same conclusions were obtained in [8] using an entirely different approach, one which is based on the use of an inner product which is not definite. Also, the criterion given below for nonnegativity is analogous to, but somewhat more involved, than the one given in [6], where the same sort of problem is considered on half line $[0, \infty)$,

THEOREM 1.3. *Assume that $\sigma > 0$, and, for $f \in F$, define*

$$D^{\sigma,\mu} f = \begin{pmatrix} f(0) - \langle f, J_0^{\sigma,\mu} \rangle \\ f(1) - \langle f, J_1^{\sigma,\mu} \rangle \end{pmatrix},$$

*where $\langle \varphi, \psi \rangle \equiv \int_0^1 \varphi(x)\psi(x)\, dx$. Then $u_f \geq 0$ if and only if $f \geq 0$ and $D^{\sigma,\mu} f = \alpha V_0^{\sigma,\mu}$ for some $\alpha \geq 0$. Moreover, if $F^{\sigma,\mu}$ denotes the subspace of $f \in F$ with $D^{\sigma,\mu} f = 0$, then $F^{\sigma,\mu}$ is invariant under $\{Q_t : t \geq 0\}$ and the restriction $\{Q_t \upharpoonright F^{\sigma,\mu} : t \geq 0\}$ is a Markov semigroup which is conservative (i.e., $Q_t \mathbf{1} = \mathbf{1}$) if* and only if $\sigma \geq \mu \coth \mu$. Finally, if $f \in F$ and $D^{\sigma,\mu} f = a_0 V_0^{\sigma,\mu} + a_1 V_1^{\sigma,\mu}$, then, uniformly for $x \in [0,1]$

(1.3) $$a_1 \neq 0 \implies \lim_{t \to \infty} e^{t\lambda_1^{\sigma,\mu}} u_f(t,x) = a_1 g_1^{\sigma,\mu}(x)$$

and

(1.4) $$a_1 = 0 \neq a_0 \implies \begin{cases} \lim_{t \to \infty} e^{t\lambda_0^{\sigma,\mu}} u_f(t,x) = a_0 g_0^{\sigma,\mu}(x), & \text{if } \sigma > \mu \coth \mu, \\ \lim_{t \to \infty} t^{-1} u_f(t,x) = a_0 g_0^{\sigma,\mu}(x), & \text{if } \sigma = \mu \coth \mu, \\ \lim_{t \to \infty} u_f(t,x) = a_0 g_0^{\sigma,\mu}(x), & \text{if } \sigma < \mu \coth \mu, \end{cases}$$

where $g_1^{\sigma,\mu}$ takes both strictly positive and strictly negative values whereas $g_0^{\sigma,\mu}$ is always strictly positive and is constant when $\sigma \leq \mu \coth \mu$. [Explicit expressions are given for $g_k^{\sigma,\mu}, k \in \{0,1\}$, in (3.1) below.]

REMARK. It should be mentioned that the Harnack principle discussed in Section 5 of [7] transfers immediately to the setting here. Namely, if $u$ is a nonnegative solution to $\dot{u} = \frac{1}{2} u'' + \mu u'$ in a region of the form $[T_1, T_2] \times [0, R]$ and $\dot{u}(t,0) = -\sigma u'(t,0)$ for $t \in [T_1, T_2]$, then, for each $T_1 < t_1 < t_2 < T_2$ and $0 < r < R$, there is a constant $C < \infty$ such that $u(s,x) \leq C u(t,y)$ for all $(s,x), (t,y) \in [t_1, t_2] \times [0, r]$, and an analogous result holds when the region is of the form $[T_1, T_2] \times [R, 1]$. The surprising aspect of this Harnack principle is that, because of the boundary condition, one can control $u(s,x)$ in terms of $u(t,y)$ even when $s \geq t$, whereas usual Harnack principles for nonnegative solutions to parabolic equations give control only when $s < t$.

1.3. *The basic probabilistic model.* The necessary stochastic calculus may be found, for example, in [2] or [5]. In particular, the second of these also contains the relevant "Markovian" results.

The probabilistic model associated with our boundary value problem can be described as follows. First, let $X$ be Brownian motion with drift $\mu$ and reflection at the boundary $\{0,1\}$. That is, if $B$ a standard Brownian motion, then one description of $X$ is as the solution to the Skorohod stochastic integral equation

$$0 \leq X_t = X_0 + B_t + \mu t + (L_0)_t - (L_1)_t \leq 1,$$

where $L_0$ and $L_1$ are the "local times" of $X$ at 0 and 1, respectively. In particular, for $k \in \{0,1\}$, $t \rightsquigarrow (L_k)_t$ is nonincreasing and increases only on $\{t : X_t = k\}$. Next, set

(1.5)
$$\Phi_t \equiv t - \sigma^{-1}(L_0)_t - \sigma^{-1}(L_1)_t,$$
$$\zeta_t \equiv \inf\{\tau > 0 : \Phi_\tau > t\} \quad \text{and} \quad Y_t \equiv X(\zeta_t).$$

When $\sigma = 0$, the interpretation of $\zeta_t$ is that it is equal $t \wedge \inf\{\tau \geq 0 : X_\tau \in \{0,1\}\}$, and so $Y$ is absorbed at the first time it leaves $(0,1)$. When $\sigma < 0$, $Y$ is Brownian motion in $(0,1)$ with drift $\mu$ and a "sticky" (i.e., it spends positive time) reflection at $\{0,1\}$. When $\sigma > 0$, $\zeta_t$ may be infinite, in which case we send $Y_t$ to a "graveyard" $\partial$ (i.e., an absorbing state outside of $[0,1]$).

The connection between (1.1) and these processes is that, for each $f \in F$ and $T \geq 0$, an application of standard Itô calculus shows that (note that $X_0 \in \{0,1\}$ and $\sigma > 0 \Longrightarrow \zeta_0 > 0$ a.s.)

(1.6) $\quad u_f(T - \Phi_t, X_t) \in \mathbb{R}$ is a continuous local martingale in $t$.

In particular,

(1.7) $\quad u_f$ bounded and $\mathbb{P}\Big(\zeta_T = \infty \Longrightarrow \lim_{t \nearrow \zeta_T} u_f(T - \Phi_t, X_t) = 0 \Big| X_0 = x\Big) = 1$
$$\Longrightarrow u_f(T, x) = \mathbb{E}[f(Y_T), \zeta_T < \infty | X_0 = x].$$

Similarly,

(1.8) $\quad u_f \geq 0 \Longrightarrow u_f(T,x) \geq \mathbb{E}[f(Y_T), \zeta_T < \infty | X_0 = x].$

REMARK. It should be emphasized that, although the process $Y$ is a familiar, continuous diffusion when $\sigma \leq 0$, it is discontinuous when $\sigma > 0$. Indeed, when $\sigma > 0$, although $Y$ behaves just like $X$ as long as it stays away from $\{0,1\}$, upon approaching $\{0,1\}$, $Y$ either jumps back inside or gets sent to $\partial$. In particular, even though it is right-continuous and has left limits, $Y$ is *not* a Hunt process because its jump times are totally accessible.

In order to make the connection between $Y$ and the functions $J_k^{\sigma,\mu}$ in Theorem 1.2, we will need the following lemma about the behavior of $\Phi_t$ as $t \to \infty$.

LEMMA 1.1. *Assume that $\sigma > 0$ and take $\mu \coth \mu = 1$ when $\mu = 0$. Then, almost surely,*

(1.9) $$\lim_{t \to \infty} \Phi_t = \begin{cases} \infty, & \text{if } \sigma > \mu \coth \mu, \\ -\infty, & \text{if } \sigma < \mu \coth \mu, \end{cases}$$

*and*

(1.10) $$\sigma = \mu \coth \mu \Longrightarrow \limsup_{t \to \infty} \pm \Phi_t = \infty.$$

*In particular, for all $T \geq 0$, $\sigma \geq \mu \coth \mu \Longrightarrow \zeta_T < \infty$ a.s. and $\sigma < \mu \coth \mu \Longrightarrow \lim_{t \to \infty} \Phi_t = -\infty$ a.s. on $\{\zeta_T = \infty\}$.*

PROOF. Assume that $\mu \neq 0$, and set

$$\psi(x) = -\left(x + \frac{e^{-2\mu x}}{\mu(1+e^{-2\mu})}\right) \coth \mu.$$

Then, $\frac{1}{2}\psi'' + \mu\psi' = -\mu \coth \mu$ and $\psi'(0) = 1 = -\psi'(1)$, and so, by Itô's formula,

$$\begin{aligned} M_t &\equiv \int_0^t \psi'(X_\tau) \, dB_\tau \\ &= \psi(X_t) + (\mu \coth \mu)t - (L_0)_t - (L_1)_t \\ &= \psi(X_t) - (\sigma - \mu \coth \mu)t + \sigma \Phi_t. \end{aligned}$$

Since $\lim_{t\to\infty} t^{-1}|M_t| = 0$ a.s., this proves that

$$\lim_{t\to\infty} \frac{\Phi_t}{t} = 1 - \frac{\mu \coth \mu}{\sigma} \qquad \text{a.s.},$$

which completes the proof of (1.9) when $\mu \neq 0$ and $\sigma \neq \mu \coth \mu$. In addition, when $\mu \neq 0$ and $\sigma = \mu \coth \mu$, the preceding says that $\psi(X_t) + \sigma\Phi_t = M_t$, and so the desired result will follow once we check that $\limsup_{t\to\infty} \pm M_t = \infty$ a.s., which, in turn, comes down to showing that $\int_0^\infty \psi'(X_\tau)^2 \, d\tau = \infty$ a.s. But, by standard ergodic theoretic considerations,

$$\lim_{t\to\infty} \frac{1}{t} \int_0^t \psi'(X_\tau)^2 \, d\tau = \int_{(0,1)} \psi'(y)^2 \nu(dy) > 0 \qquad \text{where } \nu(dy) = \frac{2\mu e^{2\mu y}}{e^{2\mu}-1} \, dy$$

is the stationary measure for $X$. Thus, the case when $\mu \neq 0$ is complete. The case $\mu = 0$ can be handled in the same way by considering the function $\psi(x) = x(1-x)$. □

As a consequence of Lemma 1.1, we can now make the connection alluded to above.

THEOREM 1.4. *Assume that $\sigma > 0$. For all bounded, measurable $\varphi: (0,1) \longrightarrow \mathbb{R}$,*

(1.11) $\qquad \mathbb{E}[\varphi(X_{\zeta_0}), \zeta_0 < \infty | X_0 = k] = \langle \varphi, J_k^{\sigma,\mu} \rangle, \qquad k \in \{0,1\}.$

*In particular, $\mathbb{P}(\zeta_0 < \infty | X_0 = k) = \langle 1, J_k^{\sigma,\mu} \rangle$ and $J_k^{\sigma,\mu}/\langle 1, J_k^{\sigma,\mu}\rangle$ is the density for the distribution of $Y_0 = X_{\zeta_0}$ given that $X_0 = k$ and $\zeta_0 < \infty$.*

PROOF. Clearly, it suffices to treat the case when $\varphi$ is continuous as well as bounded. Given such a $\varphi$, define $f \in F$ so that $f \upharpoonright (0,1) = \varphi$ and $f(k) = \langle \varphi, J_k^{\sigma,\mu} \rangle$ for $k \in \{0,1\}$. Then, by Theorem 1.3, $u_f$ is bounded and, as $t \to \infty$, $u_f(t,x) \longrightarrow 0$ uniformly for $x \in [0,1]$ when $\sigma < \mu \coth \mu$. Hence, by Lemma 1.1 and (1.7),

$$\langle \varphi, J_k^{\sigma,\mu} \rangle = f(k) = \mathbb{E}[\varphi(X_{\zeta_0}), \zeta_0 < \infty | X_0 = k]. \qquad \square$$

**2. The Riccati equation.** In this section we will prove Theorem 1.2 and the connection between solutions to (R) and solutions to (1.1). Throughout, we assume that $\sigma > 0$.

2.1. *Uniqueness of solutions to* (R).

THEOREM 2.1. *Suppose that $J \in C^2([0,1]; \mathbb{R}^2)$ is a solution to* (R), *and define $B(J)$ accordingly, as in Theorem* 1.1 *Next, for $f \in F$, set*

$$D^J f \equiv \begin{pmatrix} f(0) - \langle f, J_0 \rangle \\ f(1) - \langle f, J_1 \rangle \end{pmatrix}.$$

*Then, for any $f \in F$, $D^J u_f(t) = e^{-tB(J)} D^J f$, and so $D^J f = 0 \Longrightarrow D^J u_f(t) = 0$ for all $t \geq 0$. In particular, if $m(J) \equiv \int_0^1 |J_0(x)|\, dx \vee \int_0^1 |J_1(x)|\, dx \leq 1$, then $D^J f = 0$ implies that $\|u_f\|_{\mathrm{u}} \leq \|f\|_{\mathrm{u}}$, and, if $m(J) < 1$, then $D^J f = 0$ implies $\|u_f(t)\|_{\mathrm{u}} \longrightarrow 0$ as $t \to \infty$. Finally, if $J \geq 0$, then for any nonnegative $f \in F$ with the property that $D^J f$ is a nonnegative eigenvector of $B(J)$, $u_f \geq 0$.*

PROOF. If $J$ is any solution to (R), then,

$$\frac{d}{dt}\langle u_f(t), J \rangle = \left\langle \frac{1}{2} u_f''(t) + \mu u_f', J \right\rangle$$

$$= \left\langle u_f(t), \frac{1}{2} J'' - \mu J \right\rangle + \frac{1}{2}(u_f'(t,1)J(1) - u_f'(t,0)J(0))$$

$$- \frac{1}{2}(u_f(t,1)J'(1) - u_f(t,0)J'(0))$$

$$+ \mu(u_f(t,1)J(1) - u_f(t,0)J(0))$$

$$= -B(J)\langle u_f(t), J \rangle + \frac{d}{dt}\begin{pmatrix} u_f(t,0) \\ u_f(t,1) \end{pmatrix} + B(J)\begin{pmatrix} u_f(t,0) \\ u_f(t,1) \end{pmatrix},$$

and so $\frac{d}{dt} D^J u_f(t) = -B(J) D^J u_f(t)$, which is equivalent to $D^J u_f(t) = e^{-tB(J)} D^J f$.

Now assume that $m(J) \leq 1$ and that $D^J f = 0$. To see that $\|u_f\|_{\mathrm{u}} \leq \|f\|_{\mathrm{u}}$, let $\varepsilon > 0$ be given and suppose that $\|u_f(t)\|_{\mathrm{u}} \geq \|f\|_{\mathrm{u}} + \varepsilon$ for some $t \geq 0$. We can then find a $T > 0$ such that $\|u_f(T)\|_{\mathrm{u}} = \|f\|_{\mathrm{u}} + \varepsilon > \|u_f(t)\|_{\mathrm{u}}$ for $0 \leq t < T$. Clearly, there exists an $x \in [0,1]$ for which $|u_f(T,x)| = \|f\| + \varepsilon$. If $x \in (0,1)$, then, by the strong maximum principle for the parabolic operator $\partial_t - \frac{1}{2}\partial_x^2 - \mu\partial_x$, $|u_f|$ must be constantly equal to $\|f\|_{\mathrm{u}} + \varepsilon$ on $(0,T) \times (0,1)$, which is obviously impossible. Thus, it remains to check that $x$ can always be chosen from $(0,1)$. To this end, simply note that if $|u_f(T,x)| < \|f\|_{\mathrm{u}} + \varepsilon$ for all $x \in (0,1)$, then, for $k \in \{0,1\}$, $|u_f(T,k)| = |\langle u_f(T), J_k \rangle| < \|f\|_{\mathrm{u}} + \varepsilon$ also.

Next assume that $m(J) < 1$ and that $D^J f = 0$. To see that $\|u_f(t)\|_u \longrightarrow 0$ as $t \to \infty$, it suffices to show that $\|u_f(1)\|_u \leq \theta \|f\|_u$ for some $\theta \in (0,1)$ which is independent of $f$. Indeed, by the semigroup property and the fact that $D^J u_f(t) = 0$ for all $t \geq 0$, one would then know that $\|u_f(t)\|_u \leq \theta^n \|f\|_u$ for $t \geq n$. To produce such a $\theta$, let $\rho$ denote that first time that the process $X$ leaves $(0,1)$. Then

$$u_f(1,x) = \mathbb{E}[f(X_1), \rho > 1 | X_0 = x] + \mathbb{E}[u_f(1-\rho, X_\rho), \rho \leq 1 | X_0 = x].$$

Because $\|u_f\|_u \leq \|f\|_u$ and $|u_f(t,k)| = |\langle u_f(t,\cdot), J_k \rangle| \leq m(J) \|f\|_u$, this leads to $\|u_f(1)\|_u \leq \theta \|f\|_u$ with $\theta = 1 - \eta(1 - m(J))$, where $\eta = \inf_{x \in [0,1]} \mathbb{P}(\rho \leq 1 | X_0 = x) > 0$.

Finally, assume that $J \geq 0$ and that $D^J f$ is a nonnegative eigenvector for $B(J)$. If $f > 0$ and $u_f$ ever becomes negative, then there exists a $T > 0$ such that $u_f(t) > 0$ for $t \in [0,T)$ and $u_f(T,x) = 0$ for some $x \in [0,1]$. Again, from the strong maximum principle, we get a contradiction if $x \in (0,1)$. At the same time, because $u_f(T,k) \geq \langle u_f(T), J_k \rangle$ for $k \in \{0,1\}$, we see that the only way that $u_f(T)$ can vanish somewhere on $[0,1]$ is if vanishes somewhere on $(0,1)$. Thus, when $f > 0$, $u_f \geq 0$. To handle the case when $f \geq 0$, define $g \in F$ so that $g = 1$ in $(0,1)$ and $g(k) = \langle \mathbf{1}, J_k \rangle$ for $k \in \{0,1\}$. Next, apply the preceding result to see that $u_f + \varepsilon u_g = u_{f+\varepsilon g} \geq 0$ for all $\varepsilon > 0$, and conclude that $u_f \geq 0$. □

COROLLARY 2.1. *Let $J$ be a solution to* R *which satisfies* (1.2). *Then*

$$\langle f, J_k \rangle = \mathbb{E}[f(X_{\zeta_0}), \zeta_0 < \infty | X_0 = k] \qquad \text{for } f \in F \text{ and } k \in \{0,1\}$$

*if either $\sigma \geq \mu \coth \mu$ and (cf. the notation in Theorem 2.1) $m(J) \leq 1$ or $\sigma < \mu \coth \mu$ and $m(J) < 1$. In particular, in each of these cases, there is at most one such $J$, that $J$ must be nonnegative, and $\langle \mathbf{1}, J_k \rangle = \mathbb{P}(\zeta_0 < \infty | X = k)$ for $k \in \{1,2\}$.*

PROOF. Given the results in Theorem 2.1, there is no difference between the proof of this result and the proof given earlier of Theorem 1.4. □

By combining Theorems 1.4 and 2.1 with (1.8), we have a proof of the first assertion in Theorem 1.4. Namely, if $u_f \geq 0$, then (1.8) says that $f(k) \geq \mathbb{E}[f(X_{\zeta_0}), \zeta_0 < \infty | X_0 = k]$ and Theorem 1.4 says that $\mathbb{E}[f(X_{\zeta_0}), \zeta_0 < \infty | X_0 = k] = \langle f, J_k^{\sigma,\mu} \rangle$. Hence, we now know that $u_f \geq 0 \implies D^{\sigma,\mu} f \geq 0$, and, by the semigroup property, this self-improves to $u_f \geq 0 \implies D^{\sigma,\mu} u_f(t) \geq 0$ for all $t \geq 0$. Now suppose (cf. Theorems 1.2 and 1.3) that $D^{\sigma,\mu} f = a_0 V_0 + a_1 V_1$. Then, by Theorem 2.1, $D^{\sigma,\mu} u_f(t) = a_0 e^{-\lambda_0^{\sigma,\mu} t} V_0 + a_1 e^{-\lambda_1^{\sigma,\mu} t} V_1$. Thus, if $a_1 \neq 0$, then the ratio of the components of $D^{\sigma,\mu} u_f(t)$ is negative for sufficiently large $t > 0$, and so $a_1 = 0$ if $u_f \geq 0$. Hence, $u_f \geq 0 \implies 0 \leq D^{\sigma,\mu} f = a_0 V_0$ and therefore that $a_0 \geq 0$.

2.2. *Existence of solution to* (R). To find solutions to (R), we will first look for solutions to

(2.1) $\quad \frac{1}{2}J'' - \mu J' + BJ = 0 \quad$ with $J(0) = \begin{pmatrix} 2\sigma \\ 0 \end{pmatrix}$ and $J(1) = \begin{pmatrix} 0 \\ 2\sigma \end{pmatrix}$

for any nonsingular matrix $B$, and we will then see how to choose $B$ so that $B = B(J)$. For this purpose, set $\Omega = \sqrt{\mu^2 - 2B}$ (because of potential problems coming from nilpotence, this assignment of $\Omega$ should be thought of as an *ansatz* which is justified, ex post facto by the fact that it works) and

(2.2) $\quad J(x) = 2\sigma e^{\mu x} \left[ \frac{\sinh(1-x)\Omega}{\sinh \Omega} \begin{pmatrix} 1 \\ 0 \end{pmatrix} + e^{-\mu} \frac{\sinh x\Omega}{\sinh \Omega} \begin{pmatrix} 0 \\ 1 \end{pmatrix} \right],$

where we take $\frac{\sinh x\omega}{\sinh \omega} \equiv x$ when $\omega = 0$. It is clear that the $J$ in (2.2) solves (2.1). In addition,

$$B(J) = \sigma \left[ \mu \begin{pmatrix} -1 & 0 \\ 0 & 1 \end{pmatrix} - \Omega \coth \Omega + \frac{\Omega}{\sinh \Omega} \begin{pmatrix} 0 & e^{\mu} \\ e^{-\mu} & 0 \end{pmatrix} \right].$$

Hence, we are looking for $B$'s such that the corresponding $\Omega$ satisfies

(2.3) $\quad \frac{\mu^2 I - \Omega^2}{2\sigma} = \mu \begin{pmatrix} -1 & 0 \\ 0 & 1 \end{pmatrix} - \Omega \coth \Omega + \frac{\Omega}{\sinh \Omega} \begin{pmatrix} 0 & e^{\mu} \\ e^{-\mu} & 0 \end{pmatrix}.$

To solve (2.3), suppose that $W = (w_0, w_1)$ is a left eigenvector of $\Omega$ with eigenvalue $\omega$. Then

$$\frac{\mu^2 - \omega^2}{2\sigma} w_0 = -(\mu + \omega \coth \omega) w_0 + \frac{e^{-\mu}\omega}{\sinh \omega} w_1,$$

$$\frac{\mu^2 - \omega^2}{2\sigma} w_1 = (\mu - \omega \coth \omega) w_1 + \frac{e^{\mu}\omega}{\sinh \omega} w_0,$$

and so

$$\frac{w_1}{w_0} = \left( \frac{\mu^2 - \omega^2}{2\sigma} + \omega \coth \omega + \mu \right) \frac{e^{\mu} \sinh \omega}{\omega},$$

$$\frac{w_0}{w_1} = \left( \frac{\mu^2 - \omega^2}{2\sigma} + \omega \coth \omega - \mu \right) \frac{e^{-\mu} \sinh \omega}{\omega}.$$

In particular, $\omega$ must be a solution to

(2.4($\pm$)) $\quad \frac{\mu^2 - \omega^2}{2\sigma} + \omega \coth \omega = \pm \sqrt{\mu^2 + \frac{\omega^2}{\sinh^2 \omega}}$

and

(2.5($\pm$))
$$\frac{w_1}{w_0} = \left( \pm \sqrt{\mu^2 + \frac{\omega^2}{\sinh^2 \omega}} + \mu \right) \frac{e^{\mu} \sinh \omega}{\omega},$$

$$\frac{w_0}{w_1} = \left( \pm \sqrt{\mu^2 + \frac{\omega^2}{\sinh^2 \omega}} - \mu \right) \frac{e^{-\mu} \sinh \omega}{\omega}.$$

LEMMA 2.1. *There is a unique $\omega \geq 0$ which solves $(2.4(-))$. Moreover, if $\omega_1$ denotes this unique solution, then $\omega_1 > |\mu|$. On the other hand, $|\mu|$ is always a solution to $(2.4(+))$, and there is a second solution $\omega \in (|\mu|, \omega_1)$ if $\sigma > \mu \coth \mu$.*

PROOF. Without loss in generality, we will assume that $\mu \geq 0$.
Clearly, $\omega \geq 0$ solves $(2.4(-))$ if and only if $g_1(\omega) = 0$, where

$$g_1(\omega) \equiv \omega^2 - 2\sigma\omega \coth \omega - 2\sigma\sqrt{\mu^2 + \frac{\omega^2}{\sinh^2 \omega}} - \mu^2.$$

Since $g_1(0) < 0$ and $\lim_{\omega \to \infty} g_1(\omega) = \infty$, it is clear that $g_1$ vanishes somewhere on $(0, \infty)$. To prove that it vanishes only once and that it can do so only in $(\mu, \infty)$, first note that

$$g_1(\omega) \geq 0 \Longrightarrow (\omega - \sigma \coth \omega)^2 \geq \sigma^2 \coth^2 \omega + 2\sigma\sqrt{\mu^2 + \frac{\omega^2}{\sinh^2 \omega}} + \mu^2,$$

which is impossible unless $\omega \geq \sigma \coth \omega$, in which case $\omega > (2\sigma \coth \omega) \vee \mu$. Furthermore, if $\omega \geq 2\sigma \coth \omega$, then

$$\frac{1}{2}g_1'(\omega) = \omega - \sigma \coth \omega - \sigma \frac{1}{\sqrt{\mu^2 + \omega^2/\sinh^2 \omega}} \frac{\omega}{\sinh^2 \omega}(1 - \omega \coth \omega)$$

$$\geq \sigma \coth \omega - \frac{\sigma}{\sinh \omega} = \frac{\sigma}{\sinh \omega}(\cosh \omega - 1) > 0.$$

Knowing that $g_1(\omega) \geq 0 \Longrightarrow g_1'(\omega) > 0$ and that $\omega > \mu$, the first part of the lemma is now proved.

Turning to the second part, set

$$g_0(\omega) \equiv \omega^2 - 2\sigma\omega \coth \omega + 2\sqrt{\mu^2 + \frac{\omega^2}{\sinh^2 \omega}} - \mu^2.$$

Then $\omega$ satisfies $(2.4(+))$ if and only if $g_0(\omega) = 0$, and clearly $g_0(\mu) = 0$. In addition, since $g_1(\omega) \geq 0 \Longrightarrow g_0(\omega) > 0$ and $g_1 \geq 0$ on $[\omega_1, \infty)$, we know that $g_0$ can vanish only on $(0, \omega_1)$. Finally, to show that it vanishes somewhere on $(\mu, \omega_1)$ if $\sigma > \mu \coth \mu$, note that, since $g_0(\omega_1) > 0$ and $g_0(\mu) = 0$, it suffices to check that $\sigma > \mu \coth \mu \Longrightarrow g_0'(\mu) < 0$. But $g_0'(\mu) = (\mu \coth \mu - \sigma) \tanh \mu$, and so this is clear. $\square$

From now on, we take $\omega_1$ as in Lemma 2.1 and $\omega_0$ to be a solution to $(2.4(+))$ which is equal to $|\mu|$ if $\sigma \leq \mu \coth \mu$ and is in $(|\mu\omega_1)$ if $\sigma > \mu \coth \mu$.

The corresponding solution $J$ to (R) is given by $2\sigma e^{\mu x}/(w_{00}w_{11} - w_{01}w_{10})$ times

$$\begin{pmatrix} e^{-\mu}w_{01}w_{11}\left(\dfrac{\sinh x\omega_0}{\omega_0} - \dfrac{\sinh x\omega_1}{\omega_1}\right) + w_{00}w_{11}\dfrac{\sinh(1-x)\omega_0}{\omega_0} - w_{01}w_{10}\dfrac{\sinh(1-x)\omega_1}{\omega_1} \\ -w_{00}w_{10}\left(\dfrac{\sinh(1-x)\omega_0}{\omega_0} - \dfrac{\sinh(1-x)\omega_1}{\omega_1}\right) - e^{-\mu}w_{01}w_{10}\dfrac{\sinh x\omega_0}{\omega_0} + e^{-\mu}w_{00}w_{11}\dfrac{\sinh x\omega_1}{\omega_1} \end{pmatrix},$$

where $W_k = (w_{k0}, w_{k1})$ is a left eigenvector of $\Omega$ with eigenvalue $\omega_k$.

REMARK. For those readers who are wondering, the reason why, when $\sigma < \mu \coth \mu$, we take $\omega_0$ to be the solution to $(2.7(+))$ which is greater than $|\mu|$ is to get a solution to (R) which satisfies (1.2).

LEMMA 2.2. *The preceding $J$ is a nonnegative solution to* (R). *In addition,* $\langle 1, J_0 \rangle = 1 = \langle 1, J_1 \rangle$ *if* $\sigma \geq \mu \coth \mu$ *and* $\langle 1, J_0 \rangle \vee \langle 1, J_1 \rangle < 1$ *if* $\sigma < \mu \coth \mu$. *The eigenvalues of $B(J)$ are* $\lambda_k = \dfrac{\mu^2 - \omega_k^2}{2}, k \in \{0, 1\}$, *and associated right eigenvectors* $V_k = \begin{pmatrix} v_{k0} \\ v_{k1} \end{pmatrix}$ *satisfy*

$$\frac{v_{k1}}{v_{k0}} = (-1)^k\left(\sqrt{\mu^2 + \left(\frac{\sinh \omega_k}{\omega_k}\right)^2} + \mu\right)\frac{e^\mu \sinh \omega_k}{\omega_k}.$$

*Hence, they can be chosen so that* $v_{00} \wedge v_{01} > 0$ *with* $v_{01} + v_{01} = 1$ *and* $v_{10} > 0 > v_{11}$ *with* $v_{10} - v_{11} = 1$.

PROOF. To check that $J$ is nonnegative, we begin by remarking that $u(y) \equiv \dfrac{\sinh y\omega_0}{\omega_0} - \dfrac{\sinh y\omega_1}{\omega_1} \geq 0$ for $y \in [0, 1]$. Indeed, $u(0) = 0 = u(1)$ and $u'' \leq \omega_1^2 u$. Hence, if $u$ achieves a strictly negative minimum, it would have to do so at some $y \in (0, 1)$, in which case we would have the contradiction $0 \leq u''(y) \leq \omega_1^2 u(y) < 0$. Because of this remark, it suffices to show that all the numbers

$$\frac{w_{00}w_{11} - w_{01}w_{10}}{w_{01}w_{11}}, \quad \frac{w_{00}w_{11} - w_{01}w_{10}}{-w_{00}w_{10}},$$
$$\frac{w_{00}w_{11} - w_{01}w_{10}}{w_{00}w_{11}} \quad \text{and} \quad \frac{w_{00}w_{11} - w_{01}w_{10}}{-w_{01}w_{10}}$$

are positive. But, using $(2.8(\pm))$, this is an elementary, if somewhat tedious, task.

Next, from $B(J) = \dfrac{\mu^2 I - \Omega^2}{2}$, the identification of the eigenvalues of $B(J)$ is clear. In addition, if $W_0$ and $W_1$ are left eigenvectors of $B(J)$, then the columns of $\begin{pmatrix} W_0 \\ W_1 \end{pmatrix}^{-1}$ are associated right eigenvectors of $B(J)$. Hence, the calculation of $\dfrac{v_{k1}}{v_{k0}}$ is a consequence of $(2.8(\pm))$.

Turning to the calculation of $\langle 1, J_k \rangle$, observe that, by integrating (R), one sees that

$$B(J)\begin{pmatrix} 1 - \langle 1, J_0 \rangle \\ 1 - \langle 1, J_1 \rangle \end{pmatrix} = \mathbf{0}.$$

Hence, if $\omega_0 > |\mu|$, and therefore $B(J)$ is nondegenerate, $1 - \langle 1, J_k \rangle = 0$ for $k \in \{0, 1\}$. On the other hand, when $\omega_0 = |\mu|$, $\binom{1-\langle 1, J_0 \rangle}{1-\langle 1, J_1 \rangle}$ must be a multiple of $V_0$. In particular, this means that either $\langle 1, J_0 \rangle$ and $\langle 1, J_1 \rangle$ are both equal 1, both strictly greater than 1, or both strictly less than 1. To determine which of these holds, note that, when $\omega_0 = |\mu|$, $\frac{w_{01}}{w_{00}} = e^{2\mu}$ and therefore that

$$\langle 1, J_0 \rangle + e^{2\mu} \langle 1, J_1(x) \rangle = 2\sigma \left[ \int_0^1 e^{x\mu} \frac{\sinh(1-x)\mu}{\sinh \mu} \, dx + e^\mu \int_0^1 e^{x\mu} \frac{\sinh x\mu}{\sinh \mu} \right]$$
$$= \frac{2\sigma e^\mu \sinh \mu}{\mu},$$

and so

$$1 - \langle 1, J_0 \rangle + e^{2\mu}(1 - \langle 1, J_1 \rangle) = 1 + e^{2\mu} - \frac{2\sigma e^\mu \sinh \mu}{\mu}$$
$$= \frac{2e^\mu \sinh \mu}{\mu}(\mu \coth \mu - \sigma).$$

Thus, $\sigma = \mu \coth \mu \Longrightarrow \langle 1, J_k \rangle = 1$ and $\sigma < \mu \coth \mu \Longrightarrow \langle 1, J_k \rangle < 1$ for $k \in \{0, 1\}$. □

**3. Growth of solutions.** In this section we will give the proof of the final part of Theorem 1.3. To this end, set

$$c_k = \frac{(-1)^k \sqrt{\mu^2 \cosh^2 \omega_k + \omega_k^2 - \mu^2} - \mu \cosh \omega_k}{\omega_k + \mu} \quad \text{for } k \in \{0, 1\},$$

and define $h_0^{\sigma,\mu}$ and $h_1^{\sigma,\mu}$ by

(3.1) $$h_0^{\sigma,\mu}(x) = \begin{cases} (e^{x\omega_0} + c_0 e^{(1-x)\omega_0})e^{-x\mu}, & \text{if } \sigma > \mu \coth \mu, \\ \frac{1}{|\mu|} + \frac{x}{\mu} + \frac{1}{2\mu^2}(1 + \tanh \mu)e^{-x2\mu}, & \\ & \text{if } \sigma = \mu \coth \mu \text{ and } \mu \neq 0, \\ 1 - x(1-x), & \text{if } \sigma = 1 \text{ and } \mu = 0, \\ 1, & \text{if } \sigma < \mu \coth \mu, \end{cases}$$
$$h_1^{\sigma,\mu}(x) = (e^{x\omega_1} + c_1 e^{(1-x)\omega_1})e^{-x\mu}.$$

If $u_k^{\sigma,\mu}$ denotes $u_{h_k^{\sigma,\mu}}$, then

$$u_0^{\sigma,\mu}(t, x) = \begin{cases} e^{-t\lambda_0^{\sigma,\mu}} h_0^{\sigma,\mu}(x), & \text{if } \sigma > \mu \coth \mu, \\ t + h_0^{\sigma,\mu}(x) & \text{if } \sigma = \mu \coth \mu, \\ 1, & \text{if } \sigma < \mu \coth \mu, \end{cases}$$

and

$$u_1^{\sigma,\mu}(t, x) = e^{-t\lambda_1^{\sigma,\mu}} h_1^{\sigma,\mu}(x).$$

In addition, because $u_0^{\sigma,\mu} \geq 0$, the first part of Theorem 1.3 says that $D^{\sigma,\mu} h_0^{\sigma,\mu}$ is a nonnegative, scalar multiple of $V_0$. At the same time, because, $u_0^{\sigma,\mu}$ is unbounded when $\sigma \geq \mu \coth \mu$ and when $\sigma < \mu \coth \mu$ it does not tend to 0 as $t \to \infty$, this scalar cannot be 0. Hence, there exists a $K_0^{\sigma,\mu} > 0$ so that $K_0^{\sigma,\mu} D^{\sigma,\mu} h_0^{\sigma,\mu} = V_0$. We next want to show that $K_1^{\sigma,\mu} \neq 0$ can be chosen so that $K_1^{\sigma,\mu} D^{\sigma,\mu} h_1^{\sigma,\mu} = V_1$. It is clear (cf. Theorem 2.1) that

$$-B^{\sigma,\mu} D^{\sigma,\mu} h_1^{\sigma,\mu} = \frac{d}{dt} D^{\sigma,\mu} u_1^{\sigma,\mu}(t) \bigg|_{t=0} = -\lambda_1^{\sigma,\mu} D^{\sigma,\mu} h_1^{\sigma,\mu}.$$

Thus $D^{\sigma,\mu} h_1^{\sigma,\mu}$ is a scalar multiple of $V_1$, and, because $u_1^{\sigma,\mu}$ is unbounded, this scalar cannot be 0. That is, $K_1^{\sigma,\mu} \neq 0$ can be chosen to make $K_1^{\sigma,\mu} D^{\sigma,\mu} g_1^{\sigma,\mu} = V_1$. Finally, $h_1^{\sigma,\mu}$ must take both strictly positive and strictly negative values. If not, $u_1^{\sigma,\mu}$ would have to take only one sign, which would lead to that contradiction that $D^{\sigma,\mu} h_1^{\sigma,\mu}$ is a multiple of $V_1$.

To complete the program, set

$$g_0^{\sigma,\mu} = \begin{cases} K_0^{\sigma,\mu} h_0^{\sigma,\mu}, & \text{if } \sigma > \mu \coth \mu, \\ K_0^{\sigma,\mu}, & \text{if } \sigma \leq \mu \coth \mu, \end{cases}$$

and $g_1^{\sigma,\mu} = K_1^{\sigma,\mu} h_1^{\sigma,\mu}$. Given $f \in F$, determine $a_0$ and $a_1$ by $D^{\sigma,\mu} f = a_0 V_0 + a_1 V_1$, and set $\tilde{f} = f - a_0 g_0^{\sigma,\mu} - a_1 g_1^{\sigma,\mu}$. Then $u_f = u_{\tilde{f}} + a_0 K_0^{\sigma,\mu} u_0^{\sigma,\mu} + a_1 K_1^{\sigma,\mu} u_1^{\sigma,\mu}$. Because $D^{\sigma,\mu} \tilde{f} = \mathbf{0}$, as $t \to \infty$, $u_{\tilde{f}}(t, \cdot)$ tends to 0 if $\sigma < \mu \coth \mu$ and, in any case, stays bounded. Clearly, the last part of Theorem 1.3 follows from these considerations.

As a consequence of the preceding, we see that $-\lambda_1^{\sigma,\mu}$ is the exact exponential rate constant governing the growth of the semigroup $\{Q_t : t \geq 0\}$. That is, there is a $C < \infty$ such that

$$\|Q_t f\|_{\mathrm{u}} \leq C e^{-t \lambda_1^{\sigma,\mu}} \|f\|_{\mathrm{u}}, \tag{3.2}$$

and there are $f$'s for which $\lim_{t \to \infty} e^{t \lambda_1^{\sigma,\mu}} \|Q_t f\|_{\mathrm{u}} > 0$.

## APPENDIX

This appendix is devoted to the proof of Theorem 1.1, and we begin by introducing a little notation. First, let $g(t,x) = (2\pi t)^{-1/2} e^{-x^2/2t}$ be the centered Gauss kernel with variance $t$, and set $G(t,x) = \sum_{k \in \mathbb{Z}} g(t, x + 2k)$. Clearly, $G(t, \cdot)$ is even and is periodic with period 2. Next, set

$$Q^0(t,x,y) = e^{\mu(y-x) - \mu^2 t/2} [G(t, y-x) - G(t, y+x)], \tag{A.1}$$

$$(t,x,y) \in (0, \infty) \times [0,1]^2.$$

As one can easily check, $Q^0$ is the fundamental solution to $\dot{u} = \frac{1}{2} u'' + \mu u'$ in $[0, \infty) \times (0, 1)$ with boundary condition 0 at $\{0, 1\}$. Equivalently, if $\tau_k$ denotes $\inf\{t \geq 0 : X_t = k\}$, then

$$\mathbb{P}(X_t \in dy \text{ and } \tau_0 \wedge \tau_1 > t | X_0 = x) = Q^0(t,x,y)\, dy.$$

Next, set

$$q_k(t,x) = (-1)^k \frac{1}{2} \frac{d}{dy} Q^0(t,x,y)\bigg|_{y=k}, \qquad k \in \{0,1\}.$$

Then, by Green's theorem, for $h_k \in C([0,\infty); \mathbb{R})$,

$$w(t,x) = \int_0^t q_0(t-\tau) h_0(\tau)\,d\tau + \int_0^t q_1(t-\tau) h_1(\tau)\,d\tau$$

is the solution to $\dot{u} = \frac{1}{2}u'' + \mu u'$ in $[0,\infty) \times (0,1)$ satisfying $\lim_{t \searrow 0} u(t,\cdot) = 0$ and $\lim_{x \to k} u(t,x) = h_k(t)$. Equivalently,

$$\mathbb{P}(\tau_1 > \tau_0 \in dt | X_0 = x) = q_0(t,x)\,dt,$$
$$\mathbb{P}(\tau_0 > \tau_1 \in dt | X_0 = x) = q_1(t,x)\,dt.$$

In particular, these lead to $q_k \geq 0$ and

$$Q^0(s+t,x,y) = \int_{(0,1)} Q^0(s,x,z) Q^0(t,z,y)\,dx,$$

(A.2) $$q_k(s+t,x) = \int_{(0,1)} Q^0(s,x,y) q_k(t,y)\,dy \qquad \text{for } k \in \{0,1\},$$

$$\int_{(0,1)} Q^0(t,x,y)\,dy + \int_0^t q_0(\tau,x)\,d\tau + \int_0^t q_1(\tau,x)\,d\tau = 1$$

and

(A.3)
$$q_0(t,x) = -e^{-\mu x - \mu^2 t/2} G'(t,x),$$
$$q_1(t,x) = -e^{\mu(1-x) - \mu^2 t/2} G'(t,1-x),$$

where the second of these comes from $G'(t,1+x) = -G'(t,-1-x) = -G'(t,1-x)$.

Clearly,

(A.4) $$0 \leq Q^0(t,x,y) \leq g(t,x-y) \leq \frac{1}{\sqrt{2\pi t}}.$$

In order to estimate $q_k(t,x)$, first note that, from (A.3), it is clear that $G'(t,x) \leq 0$. Second,

$$G'(t,x) = -\frac{x}{t} G(t,x) + \frac{2}{t} \sum_{m=1}^\infty m(g(t,2m-x) - g(t,2m+x)) \geq -\frac{x}{t} G(t,x).$$

Hence,

(A.5) $$|G'(t,x)| \leq \frac{x}{t} G(t,x) \leq C \frac{x}{t} g(t \wedge 1, x)$$

and so

(A.6)
$$0 \leq q_0(t,x) \leq C\frac{x}{t} g(t \wedge 1, x),$$
$$0 \leq q_1(t,x) \leq C\frac{1-x}{t} g(t \wedge 1, 1-x)$$

for some $C < \infty$.

In what follows, we will be using the notation
$$w_0(x) = \frac{e^{2\mu x} - e^{2\mu}}{1 - e^{2\mu}}, \qquad w_1(x) = \frac{e^{2\mu x} - 1}{e^{2\mu} - 1} \quad \text{and}$$
$$\hat{f}_k = \langle f, w_k \rangle \qquad \text{for } f \in F.$$

Note that if $u \in U$ satisfies (1.1), then, after integrating by parts, one finds that
$$\dot{\hat{u}}_0(t) = -\frac{1}{2} u'(t,0) + \frac{\mu e^{2\mu}}{e^{2\mu} - 1} u(t,1) - \frac{\mu e^{2\mu}}{e^{2\mu} - 1} u(t,0),$$
$$\dot{\hat{u}}_1(t) = \frac{1}{2} u'(t,1) - \frac{\mu}{e^{2\mu} - 1} u(t,1) + \frac{\mu}{e^{2\mu} - 1} u(t,0),$$

and therefore
$$\frac{d}{dt}\begin{pmatrix} u(t,0) \\ u(t,1) \end{pmatrix} = 2\sigma \frac{d}{dt}\begin{pmatrix} \hat{u}_0(t) \\ \hat{u}_1(t) \end{pmatrix} + A\begin{pmatrix} u(t,0) \\ u(t,1) \end{pmatrix},$$

where
$$A \equiv \frac{2\sigma\mu}{e^{2\mu} - 1}\begin{pmatrix} e^{2\mu} & -e^{2\mu} \\ -1 & 1 \end{pmatrix}.$$

Solving this, we see that
$$e^{-tA}\begin{pmatrix} u(t,0) \\ u(t,1) \end{pmatrix} - e^{-sA}\begin{pmatrix} u(s,0) \\ u(s,1) \end{pmatrix}$$
$$= 2\sigma e^{-tA}\begin{pmatrix} \hat{u}_0(t) \\ \hat{u}_1(t) \end{pmatrix} - 2\sigma e^{-sA}\begin{pmatrix} \hat{u}_0(s) \\ \hat{u}_1(s) \end{pmatrix} + 2\sigma \int_s^t e^{-\tau A} A \begin{pmatrix} \hat{u}_0(\tau) \\ \hat{u}_1(\tau) \end{pmatrix} d\tau,$$

from which it is clear that if, as $s \searrow 0$, $u(s,\cdot) \restriction (0,1)$ converges pointwise to a function $f : (0,1) \longrightarrow \mathbb{R}$, then $\lim_{s \searrow 0} u(s,k)$ exists for $k \in \{0,1\}$. Thus, the first part of Theorem 1.1 is proved, and, in addition, we know that

(A.7)
$$\begin{pmatrix} u(t,0) \\ u(t,1) \end{pmatrix} = e^{tA}\begin{pmatrix} f(0) - 2\sigma\hat{f}_0 \\ f(1) - 2\sigma\hat{f}_1 \end{pmatrix} + 2\sigma\begin{pmatrix} \hat{u}_0(t) \\ \hat{u}_1(t) \end{pmatrix}$$
$$+ 2\sigma \int_0^t e^{(t-\tau)A} A \begin{pmatrix} \hat{u}_0(\tau) \\ \hat{u}_1(\tau) \end{pmatrix} d\tau$$

if $u(t,\cdot) \longrightarrow f$ in $F$.

Because, for any $u \in U$ satisfying $\dot{u} = \frac{1}{2}u'' + \mu u'$ and, as $t \searrow 0$, $u(t, \cdot) \longrightarrow f$ pointwise on $(0, 1)$,

$$\begin{aligned} u(t, x) &= \mathbb{E}[f(X_t), \sigma_0 \wedge \sigma_1 > t | X(0) = x] \\ &\quad + \mathbb{E}[u(t - \sigma_0, 0), \sigma_0 < t \wedge \sigma_1 | X(0) = x] \\ &\quad + \mathbb{E}[u(t - \sigma_1, 0), \sigma_1 < t \wedge \sigma_0 | X(0) = x] \\ &= \int_{(0,1)} Q^0(t, x, y) f(y) \, dy + \int_0^t q_0(\tau, x) u(t - \tau, 0) \, d\tau \\ &\quad + \int_0^t q_1(\tau, x) u(t - \tau, 1) \, d\tau, \end{aligned}$$

(A.7) tells us that if $u \in U$ satisfies (1.1) and $u(t, \cdot) \longrightarrow f$ in $F$, then

(A.8) $$u(t, x) = r_f(t, x) + \int_0^t k(t - \tau, x) \begin{pmatrix} \hat{u}_0(\tau) \\ \hat{u}_1(\tau) \end{pmatrix} d\tau,$$

where

$$r_f(t, x) \equiv h_f(t, x) + \int_0^t q(t - \tau, x) e^{\tau A} \begin{pmatrix} f(0) - 2\sigma \hat{f}_0 \\ f(1) - 2\sigma \hat{f}_1 \end{pmatrix} d\tau,$$

$$k(t, x) \equiv 2\sigma q(t, x) + 2\sigma \int_0^t q(t - \tau, x) e^{\tau A} A \, d\tau$$

with

$$h_f(t, x) = \int_{(0,1)} Q^0(t, x, y) f(y) \, dy \quad \text{and} \quad q(t, x) = (q_0(t, x), q_1(t, x)).$$

Our proof of the existence and uniqueness statements in Theorem 1.1 will be based on an analysis of the integral equation (A.8). Clearly, given $f \in F$, finding a solution $u$ to (A.8) comes down to finding a $t \in [0, \infty) \longmapsto v(t) = \begin{pmatrix} v_0(t) \\ v_1(t) \end{pmatrix} \in \mathbb{R}^2$ which satisfies

(A.9) $$v(t) = \hat{r}_f(t) + \int_0^t \hat{K}(t - \tau) v(\tau) \, d\tau,$$

where

$$\hat{r}_f(t) = \begin{pmatrix} \langle r_f(t, \cdot), w_0 \rangle \\ \langle r_f(t, \cdot), w_1 \rangle \end{pmatrix} \quad \text{and} \quad \hat{K}(t) = \begin{pmatrix} \langle k(t, \cdot), w_0 \rangle \\ \langle k(t, \cdot), w_1 \rangle \end{pmatrix}.$$

Indeed, if $v$ solves (A.9) and $u$ is defined by

$$u(t, x) = r_f(t, x) + \int_0^t k(t - \tau, x) v(\tau) \, d\tau,$$

then $u$ satisfies (A.8). Conversely, if $u$ solves (A.8) and $v(t) = \begin{pmatrix} \hat{u}_0(t) \\ \hat{u}_1(t) \end{pmatrix}$, then $v$ solves (A.9). Thus, existence and uniqueness for solutions to (A.8) is equivalent to existence and uniqueness for solutions to (A.9).

To prove that, for each $f \in F$, (A.9) has precisely one solution, we use the following simple lemma.

LEMMA A.1. *Suppose that $M : (0,T] \longrightarrow \mathbb{R} \otimes \mathbb{R}$ is a continuous, $2 \times 2$ matrix-valued function with the property that $L(T) = \sup_{t \in (0,T]} t^{1/2} \|M(t)\|_{\mathrm{op}} < \infty$ and that $v^0 : (0,T] \longrightarrow \mathbb{R}^2$ is a continuous function for which $\|v^0\|_{\alpha,T} \equiv \sup_{t \in (0,T]} t^\alpha |v^0(t)| < \infty$, where $\alpha \in [0,1)$. If $\{v^n : n \geq 1\}$ is defined inductively by*

$$v^n(t) = v^0(t) + \int_0^t M(t-\tau) v^{n-1}(\tau) \, d\tau, \qquad t \in (0,T],$$

*then*

$$\sup_{\tau \in [0,T]} |v^n(\tau) - v^{n-1}(\tau)| \leq \frac{(L(T)\sqrt{\pi})^n \Gamma(1-\alpha) \|v^0\|_{\alpha,T}}{\Gamma(n/2 + 1 - \alpha)} T^{n/2 - \alpha}.$$

*In particular, $\{v^n - v^0 : n \geq 1\}$ converges uniformly on $(0,T]$ to a contiguous function which tends to $0$ as $t \searrow 0$. Finally, if $v^\infty = v^0 + \lim_{n \to \infty}(v^n - v^0)$, then $v^\infty$ is the unique $v : (0,T] \longrightarrow \mathbb{R}^2$ satisfying*

$$v(t) = v^0(t) + \int_0^t M(t-\tau) v(\tau) \, d\tau \qquad \text{with } \|v\|_{\alpha,T} < \infty.$$

*In fact, there is a $C_\alpha < \infty$ such that $\|v^\infty\|_{\alpha,T} \leq C_\alpha L(T) \|v^0\|_{\alpha,T} e^{C_\alpha L(T) T}$.*

Using the estimates in (A.5) and applying Lemma A.1 with $\alpha = 0$, we now know that, for each $f \in F$, there is precisely one solution to (A.9), which, in view of the preceding discussion, means that there is precisely one solution to (A.8). Moreover, because every solution to (1.1) with initial data $f$ is a solution to (A.8), this proves that, for each $f \in F$, the only solution to (1.1) is the corresponding unique solution to (A.8); and, for this reason, in spite of our not having shown yet that every solution to (A.8) is an admissible solution to (1.1), we will use $u_f$ to denote this solution. Note that, from the last part of Lemma A.1 and our construction,

(A.11) $$\|u_f(t, \cdot)\|_{\mathrm{u}} \leq C \|f\|_{\mathrm{u}} e^{Ct}$$

for a suitable $C < \infty$.

What remains is to show that solutions to (A.8) have sufficient regularity to be an admissible solutions to (1.1) and that their dependence on $f$ is sufficiently continuous. To this end, return to (A.9), set $v^0 = \hat{r}_f(t)$ and

$$v^n(t) = v^0(t) + \int_0^t \hat{K}(t-\tau) v^{n-1}(\tau) \, d\tau.$$

Then

$$\dot{v}^n(t) = \dot{\hat{h}}_f(t) + \hat{q}(t) \begin{pmatrix} f(0) - 2\sigma \hat{f}_0 \\ f(1) - 2\sigma \hat{f}_1 \end{pmatrix} + \int_0^t \hat{K}(t-\tau) \dot{v}^{n-1}(\tau) \, d\tau,$$

where

$$\dot{\hat{h}}_f(t) = \begin{pmatrix} \langle \dot{h}_f(t,\cdot), w_0 \rangle \\ \langle \dot{h}_f(t,\cdot), w_1 \rangle \end{pmatrix} = \int_{(0,1)} \dot{\hat{Q}}^0(t,y) f(y)\, dy$$

with

$$\hat{Q}^0(t,y) = \begin{pmatrix} \langle Q^0(t,\cdot,y), w_0 \rangle \\ \langle Q^0(t,\cdot,y), w_1 \rangle \end{pmatrix}.$$

Using integration by parts, one sees that

$$\dot{\hat{Q}}^0(t,y) = \begin{pmatrix} e^{\mu y} G'(t,y) \\ e^{\mu(y-1)} G'(t,1-y) \end{pmatrix},$$

and therefore that the estimate in (A.5) together with Lemma A.1 guarantee that $\hat{u}_f(t) = \binom{(\hat{u}_f)_0(t)}{(\hat{u}_f)_1(t)}$ is continuously differentiable on $(0,\infty)$ and that

$$|\dot{\hat{u}}_f(t)| \le Ct^{-1/2} \|f\|_{\mathrm{u}} e^{Ct} \tag{A.12}$$

for some $C < \infty$. Combining this with (A.8), it follows that $u_f$ is continuously differentiable with respect to $t \in (0,\infty)$ and that

$$\dot{u}_f(t,x) = \dot{h}_f(t,x) + k(t,x)\hat{f} + q(t,x)\begin{pmatrix} f(0) - 2\sigma\hat{f}_0 \\ f(1) - 2\sigma\hat{f}_1 \end{pmatrix} + \int_0^t k(t-\tau)\dot{\hat{u}}_f(\tau)\, d\tau.$$

Since elementary estimates show that $\sup_{t>0}|t\dot{Q}^0(t,x,y)| < \infty$, we have now shown that

$$\|\dot{u}_f(t,\cdot)\|_{\mathrm{u}} \le Ct^{-1}\|f\|_{\mathrm{u}} e^{Ct} \tag{A.13}$$

for a suitable $C < \infty$.

It is clear from (A.8) that $u_f$ is differentiable on $(0,\infty) \times (0,1)$ and that

$$u'_f(t,x) = r'_f(t,x) + \int_0^t k'(t-\tau,x)\hat{u}_f(\tau)\, d\tau \qquad \text{for } (t,x) \in (0,\infty) \times (0,1).$$

The contribution of $h_f$ to $r'_f$ poses no difficulty and can be extended without difficulty to $(0,\infty) \times [0,1]$ as a smooth function. Instead, the problems come from the appearance of integrals of the form $\int_0^t q'_k(t-\tau)\psi(\tau)\, d\tau$ as $x \to k$. To handle such terms, we use (A.3) to write

$$q'_k(t,x) = -\mu q_k(t,x) + (-1)^k e^{\mu(k-x) - \mu^2 t/2} G''(t,k-x)$$
$$= -\mu q_k(t,x) + (-1)^{1-k} 2 e^{\mu(k-x) - \mu^2 t/2} \dot{G}(t,k-x).$$

The first term causes no problems. As for the second, we can integrate by parts to see that

$$\int_0^t \dot{G}(t-\tau,x)\psi(\tau)\, d\tau = G(t,x)\psi(0) + \int_0^t G(t-\tau,x)\dot{\psi}(\tau)\, d\tau.$$

Hence, by (A.12), the preceding expression for $u'_f(t,x)$ on $(0,\infty) \times (0,1)$ admits a continuous extension to $(0,\infty) \times [0,1]$. In addition, one can easily check from our earlier estimates, especially (A.12), that

$$\|u'_f(t,\cdot)\|_{\mathrm{u}} \leq Ct^{-1/2}\|f\|_{\mathrm{u}} e^{Ct} \tag{A.14}$$

for an appropriate $C < \infty$. Finally, because $u_f$ is smooth and satisfies $\dot u_f = \frac{1}{2}u'' + \mu u'$ on $(0,\infty) \times (0,1)$, we now see that $u''$ extends as a continuous function on $(0,\infty) \times [0,1]$ satisfying

$$\|u''(t,\cdot)\|_{\mathrm{u}} \leq Ct^{-1}\|f\|_{\mathrm{u}} e^{Ct} \tag{A.15}$$

for some $C < \infty$.

In view of the preceding, all that we have to do is check that $\dot u_f(t,k) = (-1)^{1-k}\sigma u'_f(t,k)$. To this end, observe that (A.8) is designed so that its solutions will satisfy

$$\begin{pmatrix} \dot u(t,0) \\ \dot u(t,1) \end{pmatrix} = 2\sigma \begin{pmatrix} \dot{\hat u}_0(t) \\ \dot{\hat u}_1(t) \end{pmatrix} + A \begin{pmatrix} u(t,0) \\ u(t,1) \end{pmatrix}$$

and that, because $\dot u = \frac{1}{2}u'' + \mu u'$,

$$2\sigma \begin{pmatrix} \dot{\hat u}_0(t) \\ \dot{\hat u}_1(t) \end{pmatrix} = \sigma \begin{pmatrix} -u'(t,0) \\ u'(t,1) \end{pmatrix} - A \begin{pmatrix} u(t,0) \\ u(t,1) \end{pmatrix}.$$

DEPARTMENT OF MATHEMATICS
MASSACHUSETTS INSTITUTE OF TECHNOLOGY
CAMBRIDGE, MASSACHUSETTS 02139
USA
E-MAIL: pang@math.mit.edu
dws@math.mit.edu